\documentclass{notices} 

\usepackage{amsfonts}
\usepackage{amsmath}
\usepackage{amscd}

\usepackage{amssymb}
\usepackage{xcolor}
\usepackage{tikz, tikz-cd}

\usepackage{savesym}
\savesymbol{udot}
\usepackage{mathabx} 

\usepackage{hyperref}

\numberwithin{equation}{section}

\newcommand\bR{\mathbb{R}}

\newcommand\bZ{\mathbb{Z}}

\newcommand\cC{\mathcal{C}}

\newcommand{\mfS}{\mathfrak{S}}

\renewcommand{\bold}[1]{\smallskip \noindent {\bf #1 }\nopagebreak}

\begin{document}

\title{What is... hierarchical hyperbolicity?}

\author{Alex Wright}

\date{}
\maketitle



\begin{figure}[h]
    \centering
    \includegraphics[width=0.85\linewidth]{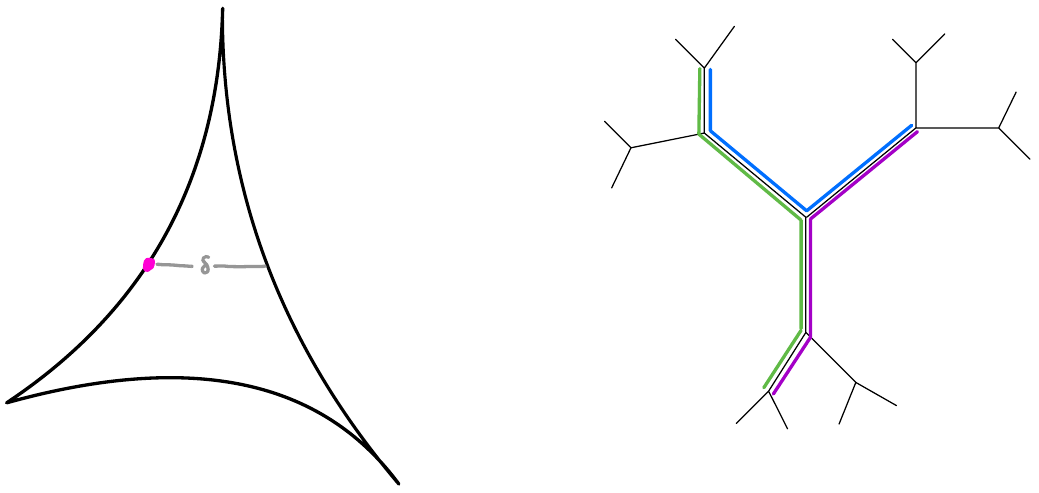}
    \caption{The definition of $\delta$-thin (left), and a triangle in a tree (right). }
    \label{F:Hyp}
\end{figure}

\bold{History.} Although curvature is defined locally, negative curvature  has implications that are ``large scale'' or ``coarse'' in that they do not reference local structure. In 1987, Gromov famously isolated one such implication from which many others follow. He defined a (Gromov) hyperbolic space to be a geodesic metric space for which there exists a $\delta\geq 0$ such that all triangles are $\delta$-thin. A triangle is by definition composed of three geodesics, and it is $\delta$-thin if each point on an edge is within $\delta$ of a point on one of the other two edges. The plane $\bR^2$ is not  hyperbolic, whereas trees have $0$-thin triangles and are the most basic examples. See Figure \ref{F:Hyp}.

A (Gromov) hyperbolic group is a group with a suitable action on a hyperbolic space. The powerful theory of hyperbolic groups unifies the algebra and geometry of free groups, fundamental groups of closed negatively curved manifolds, and many others. However, hyperbolic groups cannot have $\bZ^2$ subgroups.

Because of this, mapping class groups are not hyperbolic. It was a breakthrough when, in 1998, Masur and Minsky introduced new and powerful tools that reveal the hyperbolicity governing much of their structure \cite{MMI, MMII}. Many others added to this,  and the resulting ``Masur-Minsky machinery'' became an engine of progress in the field.

Behrstock, Hagen, and Sisto  later  gave axioms sufficient to power a type of generalized Masur-Minsky machinery, thus defining hierarchical hyperbolic groups (HHGs) and spaces (HHSs)  \cite{BehrstockHagenSistoHHSI,BehrstockHagenSistoHHSII}. 
Despite the restrictive nature of hierarchical hyperbolicity, many new examples were discovered, and the availability of a suitable framework  led to a surge of new results. Today the theory is thriving. 



\bold{Examples.} The prototypical examples of HHGs include mapping class groups and hyperbolic groups. 
The list of HHSs and HHGs now includes many cube complexes, many fundamental groups of 3-manifolds,  many Artin groups, all Teichm\"uller spaces, and the fundamental group of the moduli space of smooth cubic complex surfaces. There are many ways to produce new HHGs from old. By contrast, $SL_n(\bZ)$ is not a HHG when $n\geq 3$. 


We hope that hierarchical hyperbolicity will continue to be discovered, perhaps related to singularity theory, moduli spaces, symplectic topology, or symmetry groups of categories and cluster algebras. 

\bold{Perspectives and results.} Before turning to the “coordinate system” viewpoint, we briefly note two other perspectives for context.  
First, HHSs look hyperbolic except for product regions built from simpler HHSs. Second, combinatorial hierarchical hyperbolicity fruitfully and concisely recasts the theory in the language of simplicial complexes.

Here is a small  taste of the most easily appreciated results about HHGs, beyond what is hinted at below. 
\begin{enumerate}
\item Every element can be classified by analogy to the Nielsen-Thurston classification of mapping classes.
\item Every subgroup  either has a finite index Abelian subgroup or contains a free subgroup (strong Tits alternative).
\item Every subgroup has a finite index subgroup that, in a sense, can be block diagonalized (Omnibus Subgroup Theorem). 
\end{enumerate}

\bold{Basic setup.} Since a HHG is defined via a suitable action on a HHS, we focus primarily on spaces.


Let $(X,d_X)$ be a metric space with a collection of maps $x\mapsto x_U$, indexed by a set $\mfS$, mapping from $X$ to a collection of hyperbolic metric spaces $\cC U$.  We hope to think of this as a coordinate system, and so call $x\mapsto x_U$ a coordinate map and $x_U \in \cC U$ the $U$ coordinate of $x\in X$. Elements of the index set $\mfS$ are traditionally called ``domains.''


Distinctively, each pair of coordinate maps will  relate to each other in one of three ways, and accordingly a pair of domains can be ``nested,'' ``orthogonal,'' or ``transverse.'' 
The nesting relation is a partial order on $\mfS$, and the orthogonality and transversality relations are symmetric relations.

\bold{A key insight of hierarchical hyperbolicity.} Whereas one often tries to use the smallest number of coordinates, HHSs often use a much larger set of coordinates, whose highly structured redundancy illuminates the geometry of the space. 

The form of this highly structured redundancy would be difficult to guess without careful study of  examples. One consequence is as follows. When $U$ is transverse to or nested in $V$, there is a point $\rho^U_V\in \cC V$, and the $U$ coordinate is most useful when the $V$ coordinate is close to $\rho^U_V$. When the $V$ coordinate is not close to $\rho^U_V$,  the $U$ coordinate is, in a sense,  locally constant, and hence locally redundant.


\bold{Learning via examples.}
We pause now to examine two toy examples. Despite the fact that these two examples are atypically simple in many ways, anyone who puts in the time to understand them in detail will  find themselves surprisingly far along the road to mastering hierarchical hyperbolicity.

\begin{figure}[h]
    \centering
    \includegraphics[width=0.95\linewidth]{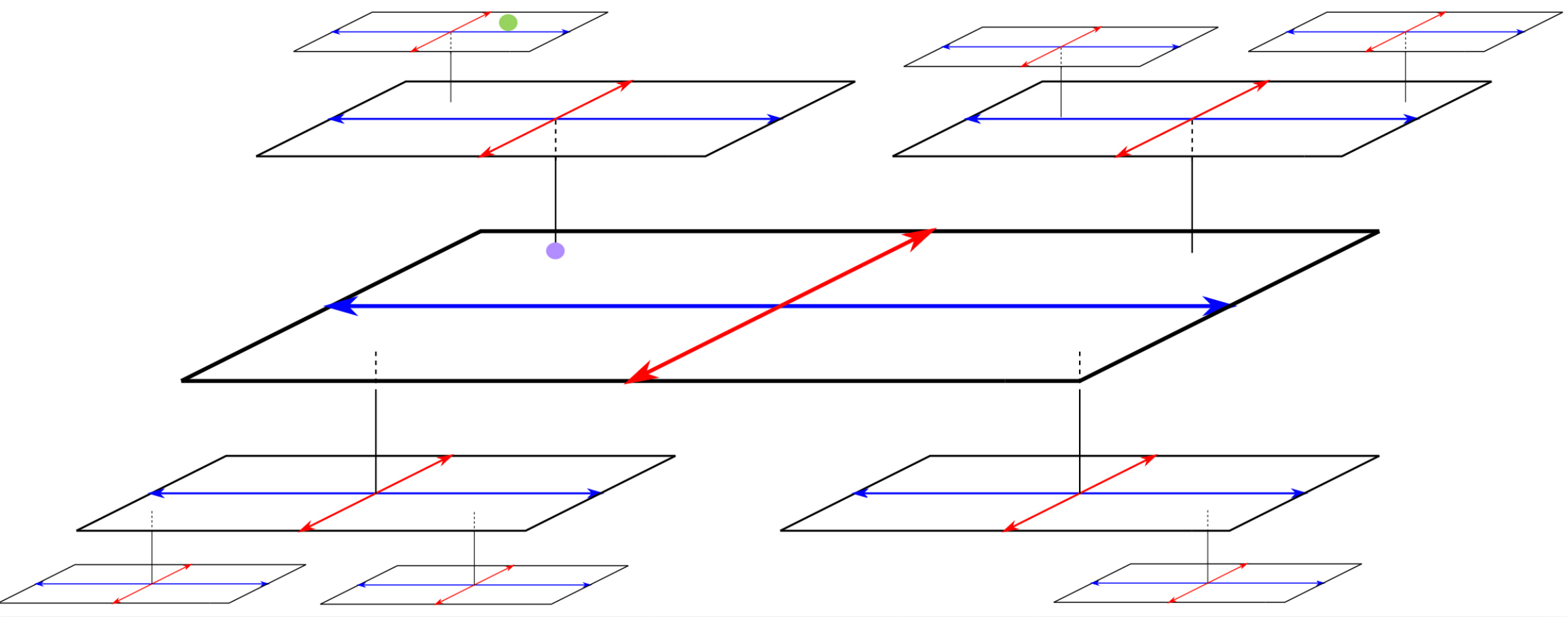}
    \caption{Jacob Russell's schematic of a tree of flats, embedded non-isometrically into $\bR^3$. The closest point projection of the green point to the centrally-drawn flat is the purple point. }
    \label{F:TreeOfFlats}
\end{figure}

\bold{First toy example.} The first toy example is illustrated in Figure \ref{F:TreeOfFlats}, and will show how tree-like and product behavior can coexist in a HHS. Let $X$ be the universal cover of the space obtained by gluing $\bR^2/\bZ^2$ to $\bR/\bZ$ at a point. This $X$ is a ``tree of flats,'' each flat being a copy of $\bR^2=\bR\times \bR$. 
We use the natural metric that is invariant under the Deck group and restricts to the $\ell^1$ metric on each flat. 

Our first domain will be denoted $T$, for ``Tree,'' and we define $\cC T$ to be the infinite-valence tree obtained by collapsing each flat to a point. The coordinate map  $X\to\cC T$ is the collapse map.

For each flat $F$ in $X$, we also define two domains $F_x, F_y$. Picking a choice of origin in $F$, define $\cC F_x=\bR$ to be the $x$-axis of $F$ and $\cC F_y=\bR$ to be the $y$-axis. (In coordinate-free terms, $\cC F_x$ is the space of vertical lines in $F$, and similarly for $\cC F_y.$)
The coordinate map to $\cC F_x$ is the $x$ coordinate of the closest point projection to the flat $F$, and similarly for $\cC F_y$. 

The $F_x$ and $F_y$ coordinates recover the information lost when collapsing the flat $F$. These coordinates are locally constant and hence redundant off $F$. It is important to record the point of $\cC T$ that is the image of $F$. We give this point two names, $\rho^{F_x}_T$ and $\rho^{F_y}_T$.

Each pair $F_x, F_y$ is orthogonal, and both $F_x, F_y$  are nested in $T$. If $E$ is another flat,  then both $F_x, F_y$ are transverse to both $E_x, E_y$. 

The region where the $E_x$ and $E_y$ coordinates can easily change, namely the flat $E$, maps via closest point projection to a single point  of $F$, again revealing  redundancy. We denote the $x$ coordinate of this point by both $\rho^{E_x}_{F_x}$ and $\rho^{E_y}_{F_x}$  and,  the $y$ coordinate by both $\rho^{E_x}_{F_y}$ and $\rho^{E_y}_{F_y}$.


The coordinates associated to $T$ and all the $F_x, F_y$ qualify as a nice system of coordinates on $X$ in the following sense:  we have a distance formula $$d_X(x,y) = \sum_{U\in \mfS} d_{\cC U}(x_U, y_U),$$ and, although $X \to \prod_{U\in \mfS} \cC U$ is not surjective, its image  admits a simple description. 

\begin{figure}[h!]
    \centering
    \includegraphics[width=0.9\linewidth]{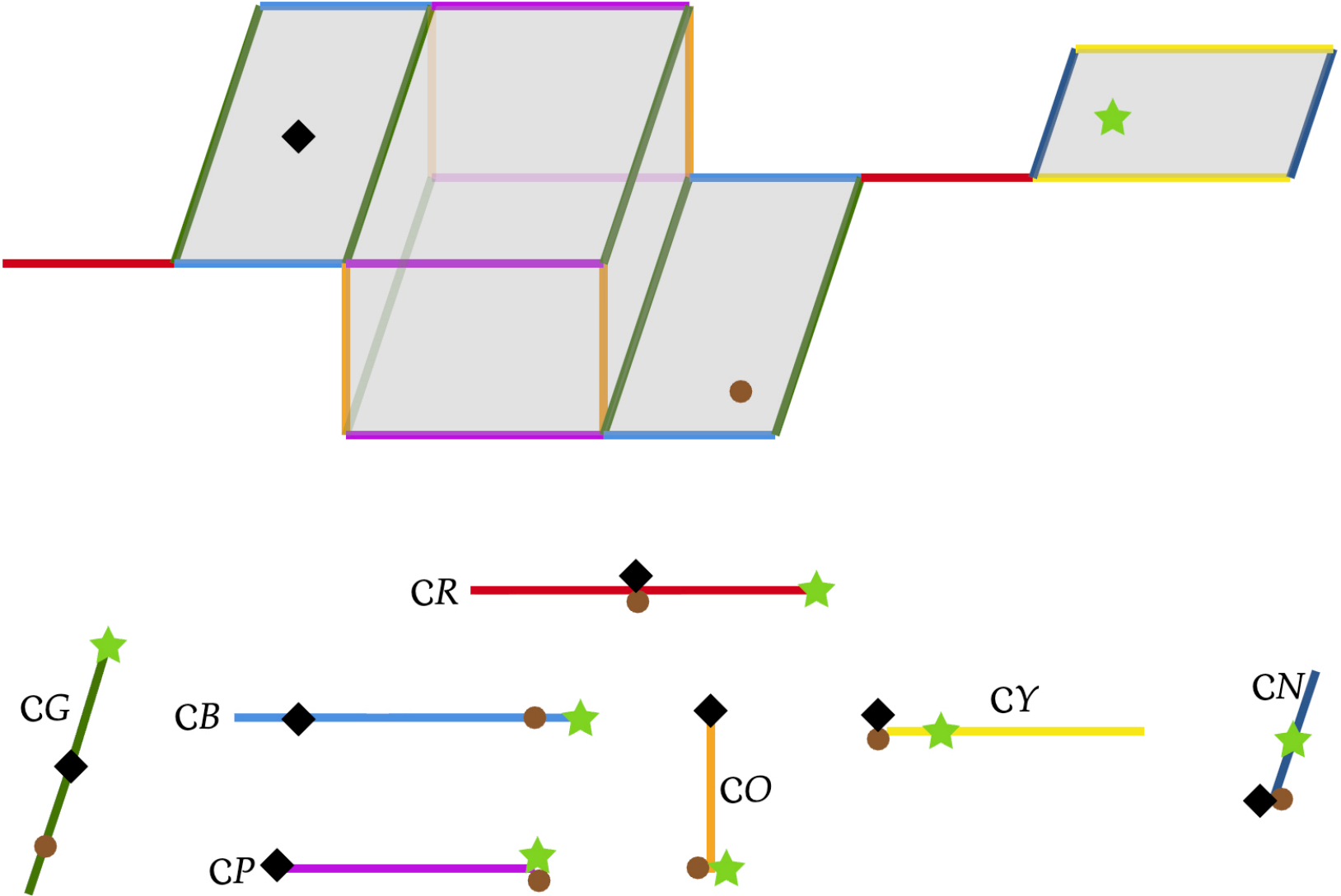}
    \caption{The second toy example (top) and the metric spaces for each domain (bottom), decorated with three arbitrary points and their coordinates.}
    \label{F:Cubical}
\end{figure}

\bold{Second toy example.} The second toy example $X$ is illustrated in Figure \ref{F:Cubical} (top), and will show how the non-tree-like parts of an HHS need not be isolated from each other. This example is a sort of cube complex obtained by gluing together products of intervals, and is significant because any two points in any HHS can be joined with a natural ``hull'' that is coarsely similar  to eamples of this form.

We define a set of coordinates indexed by $\mfS= \{R, G, B, O, P, Y, N\}$ (short for Red, Green, Blue, Orange, Purple, Yellow, and Navy). For each $U\in \mfS$ the coordinate $x_U$ takes values in an interval $\cC U$, as illustrated in the figure. 

The domain $R$ is nest-maximal, reflecting the fact that the regions where other coordinates are not redundant each get crushed to a point of $\cC R$. Similarly, $O$ and $P$ are nested in $B$. 

\begin{figure}[h!]
    \centering
    \includegraphics[width=0.9\linewidth]{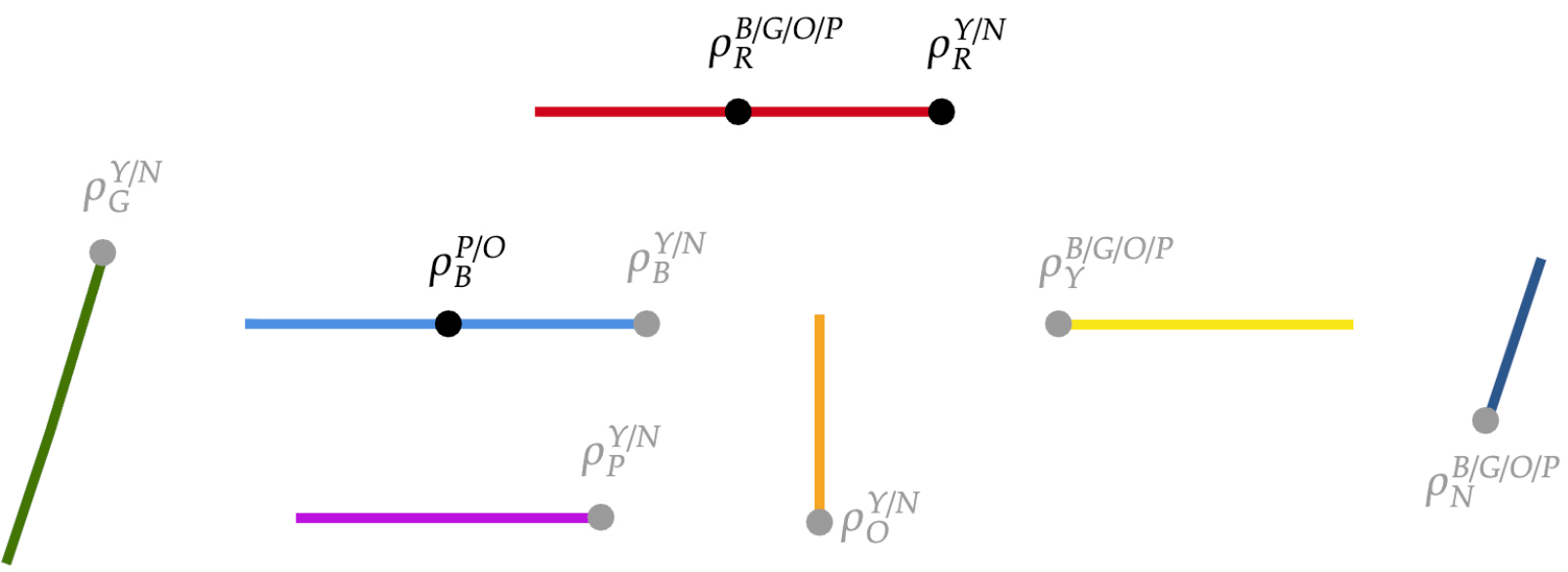}
    \caption{The $\rho$ points.}
    \label{F:CubicalRhos}
\end{figure}

There are a number of pairs, such as $\{G,B\}$ and $\{O,P\}$, for which we see a rectangle with edges of the corresponding colors. These pairs are orthogonal. Pairs of domains not related by nesting or orthogonality are  transverse. When $V$ is transverse to or nested in $U$, the point $\rho^V_U\in \cC U$ is shown in Figure \ref{F:CubicalRhos}.

\bold{Beyond toy examples.} These  examples give a good idea of what a HHS is, except that the actual definition zooms out,  discarding small-scale control to focus on large-scale geometry. For example, in general the iconic distance formula $d_X(x,y) = \sum_{U\in \mfS} d_{\cC U}(x_U, y_U)$  only holds in a coarse sense. 



\bold{Beginning the definition.} We now move on to the definition of a HHS, following \cite[Section 1.3, Remark 1.5]{BehrstockHagenSistoHHSII}. 
For expository purposes, we make a few simplifications, most notably suppressing constants and omitting the less frequently used large links axiom. We encourage the reader to focus first on the shape of the definition before mastering every detail.


We start with a geodesic  metric space $X$. We add to this a set $\mfS$ of domains, and, for each $U\in \mfS$, a hyperbolic  space $\cC U$. We also have a family $(\pi_U)_{U\in \mfS}$ of  coarsely Lipschitz maps $\pi_U: X\to \cC U$, and we think of $x_U=\pi_U(x)$ as a coordinate. 

We endow $\mfS$ with a nesting partial order $\sqsubset$ of finite height with a unique maximal element, and  an anti-reflexive, symmetric orthogonality relation $\perp$. A pair of domains that is not orthogonal and not $\sqsubset$ comparable is called transverse. We require that if $V \sqsubseteq W$ and $U\perp W$, then $U \perp V$. 

We also require that if  $V \sqsubseteq W$ then the set of domains orthogonal to $V$ and nested in $W$ is either empty or there is a domain $V^\perp \sqsubsetneq W$ that plays the role of an orthogonal complement for $V$ in $W$, in that all domains nested in $W$ and orthogonal to $V$ are nested in $V^\perp$. 

Thus far, we have a potential set of coordinates and a combinatorial structure on their index set. Now we impose structure on the coordinates, ensure we have enough coordinates, and prevent redundancy beyond what is intended. 

\bold{Highly structured redundancy.} We previously indicated that, for many pairs $x_U, x_V$ of coordinates, $x_U$ is difficult to change except when $x_V$ is close to a specific value $\rho^U_V$, suggesting that in many regions $x_U$ might be thought of as redundant. We now impose this condition with very strong structural requirements depending on  the relationship between $U$ and $V$.

\begin{figure}[h!]
    \centering
    \includegraphics[width=0.5\linewidth]{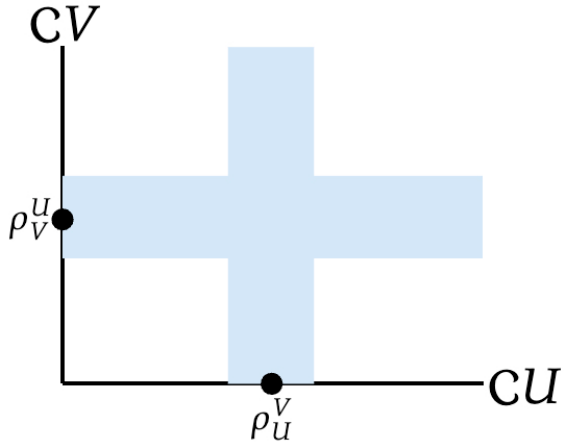}
    \caption{If $U$ and $V$ are transverse, $(x_U, x_V)$ must be in the shaded region of $\cC U \times \cC V$.}
    \label{F:Cross}
\end{figure}

If $U$ and $V$ are transverse, then there are specified points $\rho^V_U \in \cC U$ and $\rho^U_V \in \cC V$ such that for all $x\in X$ at least one of $d_{\cC U}(x_U, \rho^V_U)$ or $d_{\cC V}(x_V, \rho^U_V)$ is uniformly small. 
See Figure \ref{F:Cross}.

If $U\sqsubsetneq V$, there is a specified point $\rho^U_V\in \cC V$. One imagines $\cC U$ recaptures information lost when a subset of $X$ was crushed to the point $\rho^U_V\in \cC V$. For all $x, y \in X$ and $U\sqsubsetneq V$, if there is a geodesic in $\cC V$ from $x_V$ to $y_V$ that stays sufficiently far away from $\rho^U_V$, then $x_U$ and $y_U$ must be uniformly close.

\bold{Enough points and enough coordinates.} 
For any  pairwise orthogonal set $\{V_j\}$ and points $p_j \in \cC V_j$, there exists a point $x\in X$ such that
$d_{\cC V_j}(x_{V_j}, p_j)$ is uniformly small
for all $j$. This indicates that the $V_j$ coordinates truly do not constrain each other, and is called the partial realization axiom. 

If $d_X(x,y)$ is large, we require that there is some $U$ for which $d_{\cC U}(x_U, y_U)$ is large. 

\bold{Final  details.} We conclude with two further coherence requirements. First, if $U \sqsubsetneq V$ and both $\rho^U_W$ and $\rho^V_W$ are defined, then  
we require $d_{\cC W}(\rho^U_W, \rho^V_W)$ be uniformly small. For example, if $U \sqsubsetneq V \sqsubsetneq W$, then $\rho^U_W$ and $\rho^V_W$ should be close. 

Second, in the partial realization axiom, it should be possible to choose the point $x$ so that it lies in
the region where all the $V_j$ coordinates are easily changed.



\bold{What to read next.}  For a more in-depth introduction, see \cite{SistoWhatIs}; for recent developments, see \cite{SistoNewTools}; and for learning resources, see \cite{Website}. 

\bold{Acknowledgments.} The author gratefully acknowledges the many individuals whose comments on earlier drafts substantially improved this work.

\bibliographystyle{amsalpha}

\bibliography{bibliography}

\end{document}